\documentclass[conference]{IEEEtran}
\IEEEoverridecommandlockouts

\usepackage{cite}\usepackage{amsmath,amssymb,amsfonts,url}
\usepackage{algorithmic}
\usepackage{graphicx}
\usepackage{textcomp}
\usepackage{xcolor}
\usepackage{booktabs}
\usepackage{float}
\def\BibTeX{{\rm B\kern-.05em{\sc i\kern-.025em b}\kern-.08em
    T\kern-.1667em\lower.7ex\hbox{E}\kern-.125emX}}

\usepackage{fancyhdr}
\fancypagestyle{arxivnotice}{%
  \fancyhf{} 
  \fancyfoot[C]{\parbox{\textwidth}{\footnotesize \copyright~2026 IEEE. Personal use of this material is permitted. Permission from IEEE must be obtained for all other uses, in any current or future media, including reprinting/republishing this material for advertising or promotional purposes, creating new collective works, for resale or redistribution to servers or lists, or reuse of any copyrighted component of this work in other works.}}
}
\begin{document}

\title{Student Perceptions of Tablet-Based Teaching and Learning in TNE Undergraduate Mathematics}

\author{\IEEEauthorblockN{Dimitrios Chiotis}
\IEEEauthorblockA{\textit{School of Physical and} \\ \textit{Chemical Sciences} \\
\textit{Queen Mary University of London}\\
London, United Kingdom \\
d.chiotis@qmul.ac.uk}}

\maketitle

\thispagestyle{arxivnotice}

\begin{abstract}
Tablet-based teaching and learning has become increasingly common in university mathematics courses, offering capabilities such as digital ink, dynamic visualisation, and real-time annotation. However, there is limited experiential evidence regarding how students perceive the specific benefits and drawbacks of tablet-supported teaching and learning, particularly in Transnational Education (TNE) settings. This paper presents results from a small-scale exploratory survey of undergraduate students ($n = 19$) attending mathematics courses within an applied sciences degree programme, who anonymously  rated a set of statements regarding their experiences with tablet-based instruction on a five-point Likert scale. Initial analysis shows that students in this cohort perceived tablet teaching  most positively in supporting step-by-step understanding and active engagement during lectures  through real-time annotation and structured note-taking. In contrast, perceptions related to long-term recall, interest in mathematics, and preference over alternative teaching methods were more mixed. These findings suggest that students in this cohort perceived tablet-based teaching positively for supporting procedural understanding, while perceptions of other aspects of learning appeared more varied.
\end{abstract}


\section{Introduction}
Tablet PC devices nowadays are increasingly adopted in undergraduate mathematics courses and programmes, enabling educators to work through examples, sketch graphs, visualise functions and explain proofs. Some of the first documented findings can be traced back to Loch and Loch \& Donovan in  \cite{Loch2005,LochDonovan2006} respectively. Unlike static slides or traditional board-based teaching, tablets allow educators to write dynamically and distribute handwritten content with ease. This flexibility can assist with a clearer presentation of complex mathematical material, particularly when symbolic expressions require careful pacing, visual emphasis and step-by-step explanation, as Percival and Clayton mention in \cite{percival2015}.

Existing research tends to focus on general technology adaptations, performance outcomes, or comparisons with traditional delivery methods \cite{Haleem2022,bray2023}. There appears to be a limited  insight into what TNE students find specifically helpful or unhelpful about tablet-supported teaching and learning in mathematics courses, where symbolic reasoning and step-by-step mathematical reasoning play a central pedagogical role. 

The present study contributes to this gap by examining undergraduate students’ perceptions of tablet-based teaching and learning in one of the TNE programmes jointly delivered by Queen Mary University of London (QMUL) and Beijing University of Posts and Telecommunications (BUPT), taking place in Hainan, China. The programme follows a Flying Faculty Model (FFM) of teaching, in which academics from the home institution offer intensive teaching blocks overseas \cite{qaa2012}. Within this teaching structure, clear delivery methods and effective student engagement become particularly important during condensed teaching periods, especially in mathematics courses involving step-by-step explanations. 

Using a small-scale exploratory survey, the study identifies aspects of tablet use that students in this cohort perceived as most and least supportive. Although based on a small sample, the findings provide preliminary insights that may  inform instructional design and future research on digital technologies in mathematics education.

\section{Related work}

\subsection{Theoretical foundations and frameworks}
The integration of tablet-based learning in tertiary mathematics education is based on multiple pedagogical theories including constructivism, social learning theory, and the TPACK framework (Technological Pedagogical Content Knowledge). These frameworks emphasise the importance of aligning technology with instructional objectives and disciplinary epistemologies.

Cochrane, Narayan and Oldfield \cite{cochrane2013} present how tablets can enhance the application of social-constructivism, where students actively create together. Montrieux \cite{montrieux2018}  suggests that tablets can support active knowledge construction when the cognitive demand of navigating technology does not overwhelm learners. 

Mishra and Koehler’s TPACK framework \cite{MishraKoehler2006} suggests that effective technology integration depends on the alignment of technological, pedagogical, and content knowledge. Tablet-based teaching aligns with this perspective by combining
features such as digital ink and real-time annotation with pedagogical approaches that support mathematical explanation and interaction.



Volk et al. \cite{volk2017} evaluated the impact of cross-curricular tablet-based math activities and found that students in these environments outperformed those in traditional settings on assessments targeting higher-order learning outcomes. The findings highlight that tablet-enhanced instruction fosters skills in application, analysis, and problem-solving, aligning with the higher-levels of Bloom’s taxonomy.

Further supporting these findings, Oldfield \cite{oldfield2023} explored the role of tablets as cognitive tools in  authentic learning contexts. His doctoral study showed that mobile tablet-based learning sustained engagement over time and enabled students to apply abstract mathematical concepts to real-world problems, reinforcing the development of analytical and metacognitive skills.


 \subsection{Tablet-PCs in  teaching and learning}
Studies examining university classroom implementations have found that students perceive tablet devices as useful tools for improving interaction with learning materials and lecturers. For instance, students participating in courses that incorporated iPad-supported annotation and collaborative sharing of material reported increased participation and improved understanding of mathematical explanations, particularly when lecturers used digital handwriting to demonstrate solutions \cite{WAKEFIELD2018243}. Similarly, work examining tablet-recorded lectures in undergraduate mathematics courses found that students valued the ability to revisit lecture recordings and review problem-solving procedures at their own pace, which helped reinforce the understanding of complex mathematical concepts  \cite{yoon2011}. These findings offer the additional perspective that tablet technologies may support deeper engagement with mathematical reasoning by allowing students to revisit the procedural development of solutions at their own pace.

 More recent studies in undergraduate mathematics contexts report similar positive  perceptions, with students indicating that tablet-based teaching can enhance engagement, interaction, and participation in  learning  \cite{gouia-zarrad2015, Haddley2025}.

However, student perceptions on the incorporation of tablets in the teaching process are not entirely positive and seem to depend on how the devices are incorporated. For example, while tablet-supported activities have been associated with improved performance in procedural tasks, they may be less effective for supporting theoretical understanding \cite{WAKEFIELD2018243}, although findings might be skewed due to preferential bias.  

Overall, research findings illustrate potential benefits of tablet-based teaching for engagement and procedural understanding, while providing less detail on how students perceive different aspects of its use, particularly in TNE mathematics contexts. To explore these aspects, the present study utilises a small-scale survey to collect TNE undergraduate students' perceptions towards tablet-based teaching and learning.

\section{Methods}
\subsection{Tablet-based teaching}

Tablet-based teaching was implemented across two undergraduate mathematics courses taught by the author. Students were provided with pre-uploaded teaching materials, which included structured gapped notes, designed to be completed during the lecture. These notes were explained using live handwritten input on a stylus-enabled tablet, allowing mathematical procedures and derivations to be developed step-by-step. Students were encouraged to follow the process by completing their own notes on their devices.

This instructional approach aligns with constructivist perspectives on learning, where students actively construct understanding through engagement with tasks rather than passively receiving information \cite{cochrane2013,montrieux2018}. The use of live handwritten explanations also reflects principles of social constructivism, as the lecturer and students collaboratively develop mathematical reasoning during the lecture. In addition, presenting solutions incrementally may help manage cognitive load by presenting complex symbolic information into sequential steps. From a TPACK perspective, the approach integrates technology, pedagogy, and mathematical content in an organised way.

\subsection{Participants}

Nineteen undergraduate students  
of the TNE partnership between QMUL and BUPT in Hainan 
participated in the survey. Participation was voluntary and responses were collected anonymously. 

 Of the $19$ students who participated in the survey, $15$ of them reported that they \textit{always}  utilise Tablet PCs for learning purposes, and $4$ of them indicated that they \textit{often} use such devices. These responses suggest that participants were already familiar with tablet-based learning technologies, including accessing, downloading and annotating digital learning materials. 

\subsection{Instrument}

A Likert-scale questionnaire was developed to capture students’ perceptions of tablet-based teaching. Items were measured on a five-point scale (1 = strongly disagree to 5 = strongly agree). The questionnaire included statements relating to conceptual understanding, engagement, attention, and confidence in learning mathematics. The items were informed by themes identified in the literature on tablet-based teaching and learning in mathematics education.
\subsection{Procedure}

Students completed the survey at the end of a teaching period. Survey items focused exclusively on tablet-related teaching features and experiences. Data were analysed descriptively using item-level means, standard deviations and grand means.

\subsection{Data Analysis}

Following data cleaning, Likert-scale responses were coded numerically (1–5). Descriptive statistics were computed, including item-level and theme-level means for thematic groupings (cognition, engagement, and preferences). Given the small sample size, the analysis is limited to descriptive statistics and does not include inferential methods.

\section{Results}



\subsection{Cognitive responses}

Overall, cognitive responses yielded a theme-level mean of $3.87$, as shown in Table \ref{tab:cognitive_outcomes}. The highest rated statement was that completing gapped notes on a tablet improved students’ understanding of mathematics content (Mean $=4.00$, Standard Deviation $=1.05$). Responses relating to the lecturer’s use of the tablet showed relatively consistent agreement (Standard Deviation $=0.98$), suggesting that the manner in which the technology is used may influence student perceptions.

The lowest rated cognitive item concerned remembering material after tablet-based lectures (Mean $=3.74$, Standard Deviation $=0.99$). While this difference is modest, it may indicate that the students in this cohort perceived tablet-based teaching as supporting immediate understanding more than longer-term recall. This aligns with work in the literature suggesting that tablet-based approaches can be particularly effective for procedural engagement, while their impact on other aspects of learning may depend on how they are implemented \cite{WAKEFIELD2018243}.

\begin{table}[hbtp]
\centering
\caption{Cognition: Mean and Standard Deviation}
\label{tab:cognitive_outcomes}
\begin{tabular}{p{6.5cm} c c}
\toprule
\textbf{Statement} & \textbf{Mean} & \textbf{SD} \\
\midrule
Tablet-based teaching helps me understand mathematical concepts. & 3.95 & 1.03 \\
Tablet activities help me apply course content to solve math problems. & 3.95 & 1.03 \\
Completing gapped notes improves my understanding of the mathematics content. & 4.00 & 1.05 \\
After a tablet-based lecture, I can remember most of the taught material. & 3.74 & 0.99 \\
I can easily remember math methods written on a tablet in real time. & 3.79 & 1.03 \\
The lecturer uses the tablet effectively to explain mathematical concepts. & 3.79 & 0.98 \\
Overall, tablet-based instruction enhances my mathematics learning. & 3.84 & 1.12 \\
\hline \\
\textbf{Theme-level Mean} &3.87 &\\ 
\bottomrule
\end{tabular}
\end{table}

\subsection{Engagement}

Engagement‑related statements are presented in Table \ref{tab:affective_outcomes}. The highest rated response was about following the lecture better when annotating  gapped notes in real-time on a tablet (Mean = 4.44, Standard Deviation = 0.73). Students also reported improved confidence in learning mathematics with tablet‑based teaching (Mean = 4.16, Standard Deviation = 1.01). These findings are consistent with earlier work in \cite{Loch2005} and \cite{LochDonovan2006}, who noted that real‑time handwriting on tablets may help students stay engaged with a mathematics lecture. In \cite{GalliganLoch2010}, it is similarly highlighted that the ability to annotate alongside the lecturer enhances both attention and understanding.

The lowest rated engagement items were tablet‑based teaching increasing interest in learning mathematics (Mean = 3.74, Standard Deviation = 1.15) and tablet‑based mathematics teaching being more engaging than non‑tablet alternatives (Mean = 3.74, Standard Deviation = 1.15). The relatively high standard deviations here suggest that students’ experiences were mixed. This affirms Wakefield et al. \cite{WAKEFIELD2018243}, who found that while tablets can enhance engagement, they do not necessarily make a subject more interesting.   A. \& J. Haddley \cite{Haddley2025} also observed that even students who strongly prefer blackboards still valued multi‑modal delivery, pointing to the importance of how technology is integrated rather than the device itself.

\begin{table}[htbp]
\centering
\caption{Engagement: Mean and Standard Deviation}
\label{tab:affective_outcomes}
\begin{tabular}{p{6.5cm} c c}
\toprule
\textbf{Statement} & \textbf{Mean} & \textbf{SD} \\
\midrule
Tablet-based teaching improved my confidence in writing mathematics. & 3.95 & 1.08 \\
Tablet-based teaching improved my confidence in learning mathematics. & 4.16 & 1.01 \\
Gapped notes help me stay concentrated and follow the lecture. & 3.89 & 1.15 \\
I can follow the lecture better when I fill in gapped notes. & 4.44 & 0.73 \\
The lecturer uses the tablet in an engaging and interactive way. & 4.00 & 0.71 \\
The lecturer's tablet use makes mathematics lectures feel more personal and dynamic. & 3.89 & 0.99 \\
The lecturer's enthusiasm for technology positively influences my engagement in maths. & 3.84 & 1.01 \\
I pay greater attention when real-time handwriting on a tablet is used. & 3.84 & 1.12 \\
Tablet-based teaching increases my interest in learning mathematics. & 3.74 & 1.15 \\
Tablet-based mathematics teaching is more engaging than non-tablet based. & 3.74 & 1.15 \\
\hline \\
\textbf{Theme-level Mean} &3.95 &\\ 
\bottomrule
\end{tabular}
\end{table}

\subsection{Preferences}
Preference statements are given in Table \ref{tab:comparative_outcomes}. Students showed moderate agreement with wanting more tablet‑based courses (Mean = 3.74, Standard Deviation = 0.99) and with tablets being more engaging than non‑tablet alternatives (Mean = 3.74, Standard Deviation = 1.15). The moderate means and higher variability align with \cite{Haddley2025}, where it was found that students often prefer multimodal delivery over a single mode, and with \cite{WAKEFIELD2018243}, where it was noted that tablets are valued for specific tasks but do not automatically replace other approaches. 
\begin{table}[htbp]
\centering
\caption{Preferences: Mean and Standard Deviation}
\label{tab:comparative_outcomes}
\begin{tabular}{p{6.5cm} c c}
\toprule
\textbf{Statement} & \textbf{Mean} & \textbf{SD} \\
\midrule
Tablet-based mathematics teaching is more engaging than non-tablet based. & 3.74 & 1.15 \\
Overall, tablet-based instruction enhances my mathematics learning. & 3.84 & 1.12 \\
I would prefer more mathematics courses taught using tablets. & 3.74 & 0.99 \\
\bottomrule
\end{tabular}
\end{table}

\subsection{Most and least helpful aspects}

Across all items, the most highly rated aspects of tablet-based teaching were related to following step-by-step explanations and actively completing gapped notes during lectures. These features were consistently associated with higher perceived understanding, engagement, and confidence in learning mathematics. 

In contrast, lower-rated items were associated with increased interest in mathematics, preference for tablet-based instruction over alternative methods, and perceived support for memory and recall. These initial results suggest that students in this cohort perceived tablet-based teaching as particularly helpful for supporting procedural understanding in real time, while their perceptions about broader learning preferences appeared to be more variable.

\section{Discussion}

The preliminary findings of this small-scale study suggest that tablet-based teaching was perceived by students in this cohort as  supportive of procedural understanding and engagement during mathematics lectures. The highest-rated items for cognition and engagement were  associated with following step-by-step explanations and actively completing gapped notes. This aligns with prior work \cite{Loch2005,LochDonovan2006,GalliganLoch2010}, which emphasises the value of real-time handwritten exposition in supporting mathematical reasoning and maintaining student attention and engagement.

At the same time, lower-rated items were associated with  interest in mathematics and preference for tablet-based instruction over alternative methods. These findings are consistent with studies such as \cite{WAKEFIELD2018243,Haddley2025}, which suggest that while tablet technologies can enhance certain aspects of engagement, they do not necessarily replace traditional approaches or universally increase student motivation. 

The early findings of this paper highlight that the perceived effectiveness of tablet-based teaching depends not only on the technology itself but also on how it is implemented. Items relating to the lecturer’s use of the tablet showed relatively consistent agreement,  suggesting that pedagogical practice may play a central role in shaping student experience. This aligns with TPACK-informed perspectives, where the integration of technology, pedagogy, and content is considered important for effective teaching.

In the context of transnational education, these findings provide provisional insights into how students experienced technology-supported instruction across potentially diverse educational expectations. The preference for multimodal delivery observed within this cohort may reflect the need for flexible teaching approaches that accommodate varied learning backgrounds and prior experiences. However, given the small sample size and exploratory nature of the study, these interpretations should be considered indicative rather than generalisable.
\section{Conclusion}

This study examined the perceptions of $19$ undergraduate students regarding tablet-based teaching in a transnational mathematics education context. Exploratory insights suggest that tablet use is perceived as particularly helpful for supporting step-by-step understanding and active engagement during lectures, while its influence on longer-term recall and overall teaching preferences appears more variable.

These results highlight the importance of considering how specific features of tablet-based instruction are experienced by students within this specific cohort. In TNE settings, this may be particularly relevant, as students encounter diverse teaching approaches and expectations.

 Future work could extend this study by including larger samples and combining perception data with performance-based measures to further explore the relationship between student perceptions and actual learning outcomes.

\vspace{-1ex}

\bibliography{references}

@inproceedings{Loch2005,
  author    = {Birgit Loch},
  title     = {Tablet Technology in First Year Calculus and Linear Algebra Teaching},
 booktitle =  {5th Southern Hemisphere Conference on Undergraduate Mathematics and Statistics Teaching and Learning},

  address   = {Fraser Island, Australia},
year = {2005},
  pages     = {231--237},
}

@article{LochDonovan2006,

  author  = {Birgit Loch and Dianne Donovan},

  title   = {Progressive Teaching of Mathematics with Tablet Technology},

  journal = {Journal of Instructional Science and Technology (e-JIST)},

  volume  = {9},

  number  = {2},

  year    = {2006}

}

@article{GalliganLoch2010,
    author  = {Linda Galligan and Birgit Loch and Christine McDonald and Janet A. Taylor},
  title     = {The Use of Tablet and Related Technologies in Mathematics Teaching},
 journal = {Australian Senior Mathematics Journal},
  volume  = {24},

  number  = {1},

  pages   = {38--51},

  year    = {2010}
}

@article{Haleem2022,
  author  = {Abid Haleem and Mohd Javaid and Mohd Asim Qadri and Rajiv Suman},
  title   = {Understanding the Role of Digital Technologies in Education: A Review},
  journal = {Sustainable Operations and Computers},
  volume  = {3},
  pages   = {275--285},
  year    = {2022}
}

@article{MishraKoehler2006,

  author  = {Punya Mishra and Matthew J. Koehler},

  title   = {Technological Pedagogical Content Knowledge: A Framework for Teacher Knowledge},

  journal = {Teachers College Record},

  volume  = {108},

  number  = {6},

  pages   = {1017--1054},

  year    = {2006},

  doi     = {10.1111/j.1467-9620.2006.00684.x}

}

@article{cochrane2013,
  author  = {Cochrane, Thomas and Narayan, Vickel and Oldfield, James},
  title   = {i{P}adagogy: Appropriating the i{P}ad within pedagogical contexts},
journal = {International Journal of Mobile Learning and Organisation},  volume  = {7},
  number  = {1},
  pages   = {48--65},
  year    = {2013}
}

@inproceedings{montrieux2018,
  author    = {Montrieux, H. and Schellens, T.},
  title     = {The impact of tablet devices on high school students' cognitive load and learning},
booktitle =  {Proceedings of the 12th International Technology, Education and Development Conference},  pages     = {1591--1596},
  address ={Valencia, Spain},
  year      = {2018}
}

@phdthesis{oldfield2023,
  author = {J. D. Oldfield},
  title  = {Mobile Authentic Learning: Tablet Devices as Cognitive Tools Supporting Authentic Learning Experiences},
  school = {Murdoch University},
  year   = {2023},
  type   = {Ph.{D}. dissertation}
}

@article{volk2017,
  author  = {Volk, M. and Cotič, M. and Zajc, M. and Starčič, A. I.},
  title   = {Tablet-based cross-curricular maths vs. traditional maths classroom practice for higher-order learning outcomes},
journal = {Computers and Education},  volume  = {114},
  pages   = {1--23},
  year    = {2017}
}

@article{bray2023,
  title     = {Mobile learning in mathematics education: A systematic literature review of empirical research},
  author = {D. M. Tang and others},
 journal = {Eurasia Journal of Mathematics, Science and Technology Education},
  volume    = {19},
  number    = {5},
  pages     = {em2268},
  year      = {2023}
}

@article{yoon2011,
  author  = {Susan Yoon and John Sneddon},
  title   = {Student Perceptions of Tablet-Based Lecture Recording in Undergraduate Mathematics},
 journal = {International Journal of Mathematical Education in Science and Technology},
  volume  = {42},
  number  = {4},
  pages   = {425--445},
  year    = {2011}
}

@misc{qaa2012,
  author       = {{Quality Assurance Agency for Higher Education}},
  title        = {Review of {UK} Transnational Education in {C}hina: Overview},
  year         = {2012},
 
     howpublished = {QAA Report. [Online]. Available: QAA website},}

@article{WAKEFIELD2018243,
title = {The impact of an i{P}ad-supported annotation and sharing technology on university students' learning},
journal = {Computers and Education}, volume = {122},
pages = {243-259},
year = {2018},
issn = {0360-1315},
doi = {https://doi.org/10.1016/j.compedu.2018.03.013},
url = {https://www.sciencedirect.com/science/article/pii/S0360131518300666},
author = {James Wakefield and Jessica K. Frawley and Jonathan Tyler and Laurel E. Dyson},

}

@article{Haddley2025,
  author  = {Alena Haddley and Joel Haddley},
  title   = {Student Perceptions of {T}ablet{PC} Use in Mathematics Teaching, and Student Preferences of Different Delivery Modes},
  journal = {MSOR Connections},
   note    = {vol. 24, no. 1, Nov. 2025},
  url     = {https://journals.gre.ac.uk/index.php/msor/article/view/1598}
}

@incollection{gouia-zarrad2015,

  author    = {Rim Gouia and Cindy Gunn and Della Audi},

  title     = {Using iPads in University Mathematics Classes: What Do the Students Think?},

  booktitle = {Assessing the Role of Mobile Technologies and Distance Learning in Higher Education},

  editor    = {Patricia Ord{\'o}{\~n}ez de Pablos and others},

  publisher = {IGI Global Scientific Publishing},

  pages     = {60--77},

  year      = {2015},

  doi       = {10.4018/978-1-4666-7316-8.ch003}

}

@inproceedings{percival2015,
  author    = {J. Percival and T. Claydon},
  title     = {A Study of Student and Instructor Perceptions of {Tablet PCs} in Higher Education Contexts},
 booktitle = {Higher Education in Transformation Conference},
  address   = {Dublin, Ireland},
  pages     = {250--264},
  year      = {2015}
}

\end{document}